\documentclass{article}
\usepackage{url}
\usepackage{graphicx}
\usepackage{amsmath,amsthm,amsfonts,amssymb,latexsym,enumerate}
\newtheorem{thm}{Theorem} 
\newtheorem{lem}[thm]{Lemma}

\newtheorem{conj}[thm]{Conjecture}

\newtheorem{example}[thm]{Example}

\newcommand{\Irr}{{\mathrm {Irr}}}

\newcommand{\ZZ}{{\mathbb Z}}

\newcommand{\PC}{({\mathcal P})}

\newcommand{\FF}{{\mathbb F}}

\newcommand{\SSS}{{\sf S}}

\begin{document}
\title{Machine-Assisted Proofs\\ICM 2018 Panel}
\author{James Davenport (moderator)\footnote{Department of Computer Science, University of Bath, Bath BA2 7AY, U.K., {\tt J.H.Davenport@bath.ac.uk}}{ } and Bjorn Poonen\footnote{MIT Department of Mathematics
32 Vassar St., Bldg. 2-243
Cambridge, MA 02139, USA, {\tt poonen@math.math.edu}}{ } and James Maynard\footnote{Professor of Number Theory, University of Oxford, UK}{ } and \\Harald Helfgott\footnote{Mathematisches Institut
Bunsenstra{\ss}e 3--5
D-37073 G\"ottingen, Germany;   IMJ-PRG, UMR 7586,
  58 avenue de France, B\^{a}timent S. Germain, case 7012,
  75013 Paris CEDEX 13, France, {\tt harald.helfgott@gmail.com}.}{ } and Pham Huu Tiep\footnote{Department of Mathematics,  Rutgers University, Piscataway, NJ 08854, USA 
{\tt tiep@math.rutgers.edu}}{ } and Lu\'\i{}s Cruz-Filipe\footnote{Department of Mathematics and Computer Science, University of Southern Denmark, Campusvej 55,
5230 Odense M,
Denmark, {\tt lcfilipe@gmail.com}}}
\date{7 August 2018}
\maketitle\noindent
This submission to arXiv is the report of a panel session at the 2018 International Congress of Mathematicians (Rio de Janeiro, August).  It is intended that, while \verb+v1+ is that report, this stays a living document containing the panelists', and others', reflections on the topic.
\par
This panel took place on Tuesday 7th August 2018. After the moderator had introduced the topic, the panelists presented their experiences and points of view, and then took questions from the floor.
\section{Introduction (James Davenport)}
\subsection{A (very brief, partial) history}
\begin{description}
	\item[1963]``Solvability of Groups of Odd Order'': 254 pages\footnote{``one of the longest proof to have appeared in the mathematical literature to that point.'' \cite{Gonthieretal2013a}.} \cite{FeitThompson1963}. Also Birch \& Swinnerton--Dyer published \cite{BirchSwinnertonDyer1963}, the algorithms underpinning their conjectures.
\item[1976]``Every Planar Map is Four-Colorable'': 256 pages + computation \cite{AppelHaken1976a}.
\item[1989]Revised Four-Color Theorem proof published \cite{AppelHaken1989}.
\item[1998]Hales announced proof of Kepler Conjecture.
\item[2005]Hales' proof published in an abridged form ``uncertified''\footnote{\emph{Mathematical Reviews} states ``Nobody has managed to check all the details of the proof so far, but the theoretical part seems to be correct. The whole proof is considered and assumed to be correct by most of the mathematical community.'' \url{https://mathscinet.ams.org/mathscinet-getitem?mr=2179728}.} \cite{Hales2005}.
\item[2008]Gonthier stated formal proof of Four-Color Theorem \cite{Gonthier2008}.
\item[2012]Gonthier/Th\'ery stated\footnote{``Both the size of this proof and the range of mathematics involved make formalization
a formidable task'' \cite{Gonthieretal2013a}.} formal proof of Odd Order Theorem \cite{GonthierThery2012a,Gonthieretal2013a}.
\item[2013]Helfgott published (arXiv) proof of ternary Goldbach Conjecture \cite{Helfgott2013a}.
\item[2014]Flyspeck project announced formal proof of Kepler Conjecture \cite{Hales2014a}.
\item[2015]Maynard published ``Small gaps between primes'' \cite{Maynard2015a}.
\item[2017]Flyspeck paper published \cite{Halesetal2017a}.
\end{description}
The Odd Order Theorem is important, but chiefy because it leads to the classification of finite simple groups. One might ask when this will be formally proved, and indeed I did ask Georges Gonthier this question. He answered that he worked, not so much from \cite{FeitThompson1963} itself as from \cite{Benderetal1994,Peterfalvi2000}, two substantial books which between them described much work simplifying  and clarifying the argument, and that such work had yet to be done for the full classification.
\subsection{Questions for Consideration}
What are the implications for
\begin{itemize}
\item authors
\item journals and their publishers
\item the refereeing process (we note that, although \cite{Hales2005} took seven years not to be completely refereed,  \cite{Halesetal2017a} still took three years to be refereed: ``formal'' is not the same as ``simple''.) 
\item readers
\item the storage and curation of such proofs, and, if necessary, the software necessary to run such proofs.?
\end{itemize}
These questions are not independent: the refereeing process is run by journals, and one can ask whether the journal should keep the machine-readable proof, as with \cite{Hales2005}, or whether Helfgott is right with ``available on request'', or maybe Maynard's ``at \url{www.arxiv.org}''.
\par
{\bf Acknowledgements.} I am extremely grateful to Ingrid Daubechies for her hard work convening this panel.

\section{Bjorn Poonen}
\def\green#1{{\bf #1}}
\def\magenta#1{{\it #1}}
\subsection {Kinds of machine assistance}
There are various forms of machine-assisted proof, and they need different approaches.
\begin{itemize}
\item \green{Experimental mathematics}: Humans design experiments for the computers
to carry out, in order to discover or test new conjectures.
\item \green{Human/machine collaborative proofs}: Humans reduce a proof to a large 
number of cases, or to a detailed computation, which the computer carries out.
\item \green{Formal proof verification}: Humans supply the steps of a proof
in a language such that the computer can verify that each step 
follows logically from previous steps.
\item \green{Formal proof discovery}: The computer searches for a chain
of deduction leading from provided axioms to a theorem.
\end{itemize}
\magenta{Cloud computing services} make it possible to do computations 
much larger than most people would be able to do with their 
own physical computers.

\begin{example}
My colleague Andrew Sutherland at MIT 
did a 300 core-year computation in 8 hours
using (a record) 580,000 cores of the Google Compute Engine 
for about \$20,000.
\rm 
As he says, 
\begin{quote}
``Having a computer that is 10 times as fast 
doesn't really change the way you do your research, but having a computer
that is 10,000 times as fast does.''
\end{quote}
\end{example}

\subsection {Large databases}
There are now many databases of mathematical facts, such as the Atlas of Finite Simple Groups \cite{Conwayetal1985} and the Online Encyclopedia of Integer Sequences \cite{Sloane2003,Sloane2007}. My personal tool is the LMFDB --- The L-functions and modular forms database (\url{www.lmfdb.org}).
\begin{figure}
	\caption{Modular Forms Database}
\includegraphics[scale=0.25]{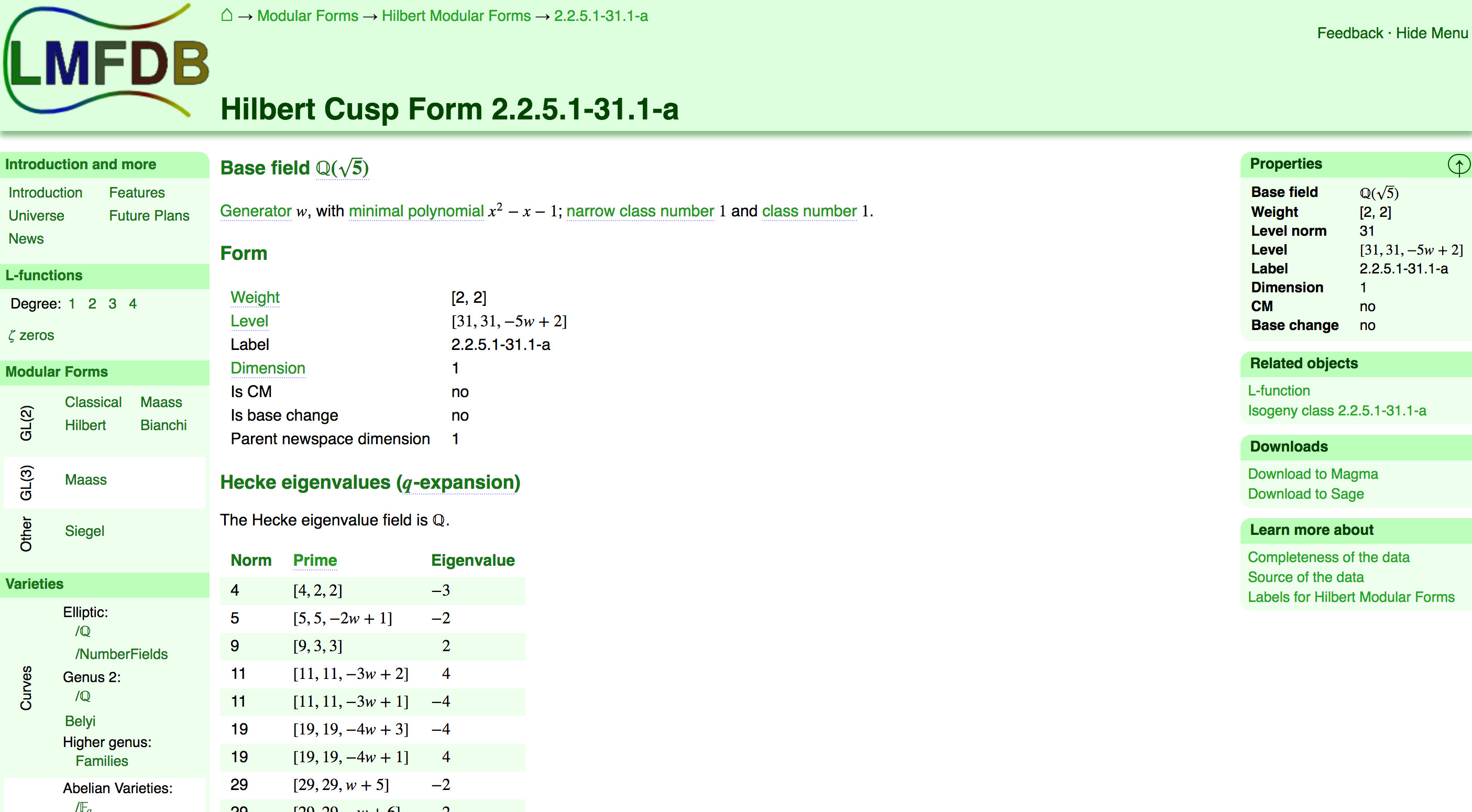}
\end{figure}

\subsection {Refereeing computations}

Almost all math journals now publish papers depending on machine computations.
\begin{itemize}
\item
Should papers involving machine-assisted proofs 
be viewed as more or less suspect than purely human ones?
\item 
What if \magenta{the referee is not qualified} to check the computations,
e.g., by redoing them?
\item 
What if the computations are \magenta{too expensive} to do more than once?
\item 
What if the computation involves \magenta{proprietary software} (e.g., MAGMA)
for which there is no direct way to check that the algorithms
do what they claim to do?
\item 
Should one insist on \green{open-source software},
or even software that has gone through a \green{peer-review process}
(as in Sage, for instance)?
\item
In what form should computational proofs be published?
\end{itemize}

\subsection {Some opinions}

\begin{itemize}
\item
The burden should be on authors 
to provide understandable and verifiable code,
just as the burden is on authors to make their human-generated proofs 
understandable and verifiable.
\item 
Ideally, code should be written in a high-level computer algebra package
whose language is close to mathematics,
to minimize the computer proficiency demands on a potential reader.
\end{itemize}

{\bf Acknowledgements.} Many thanks are due to my colleague Andrew V. Sutherland.

\section{James Maynard}
\subsection{Opportunities and challenges of use of machines}
The use of computers in mathematics is widespread and likely to increase.\\
\hfill\\
This presents several \textbf{opportunities}:
\begin{itemize}
\item Personal assistant: Guiding intuition, checking hypotheses
\item Theorem proving: Large computations, checking many cases
\item Theorem checking: Formal verification
\end{itemize}
But equally it presents several \textbf{challenges}:
\begin{itemize}
\item Are computations rigorous?
\item Are proofs with computation understandable to humans?
\item How do we find errors?
\end{itemize}

\subsection{My use of computation}
I use computation daily to guide my proofs and intuition.
Two particular kinds of computations occur in my proofs.
\par
The first is when I want to understand the spectrum of infinite dimensional operator, but of course computers can't handle these as such. Then I consider well-chosen finite-dimensional subspace and do explicit computations there.
\begin{itemize}
\item Do non-rigorous computation to guess good answers, then compute answers using exact arithmetic as in \cite{Maynard2015a}.
\item Happy trade-off between quality of numerical result and computation time.
\end{itemize}
The second is when, after doing theoretical manipulations, I need to show that some messy explicit integral is less than 1, as in \cite{Maynard2013a}.
\begin{itemize}
\item My computations are non-rigorous!
\item Very difficult to referee --- many potential sources for error.
\end{itemize}

\section{Harald Helfgott}
\subsection{The varieties of machine assistance: Main issues.}

When we say that our work on a proof is ``machine-assisted'', we may
mean any of several related things. First, we may simply refer to exploratory
work -- computations and plots that give us an idea of which statements are
true. We may also speak of computations that are truly part of a proof --
dealing with case-work, say, or numerical computations, which of course ought
to be rigorous. Lastly, we may speak of {\em formal proofs}: their production
and verification could in principle be done by hand, but part of their point
is that they can be verified mechanically; moreover,
as far as the working mathematician is concerned, they were essentially a
useful fiction before computer assistance in their creation became possible.

My present concerns center on the second kind of machine assistance just listed,
namely, the use of computations -- numerical or symbolic -- as part of a proof. 
The first issue to discuss is what is meant by rigorous computation. (That
would be less of an issue for exploratory work, where all one wants is results
that are sufficiently precise and almost certainly correct.) Another issue
is that of error. Are computer errors at all likely? How do we reduce the
possibility of programming errors, or, more generally, errors in machine-human
interaction? For that matter, can computers help us check for errors in proofs
altogether, whether the process of proof was computer-assisted or not?

This last question does lead naturally to the matter of formal proofs. That
matter gives rise to another sort of concern: what is the meaning
and purpose of a proof, for us, as human beings practicing mathematics? I will
not do more than touch upon
that domain. What can be discussed more factually is the
extent to which formal proofs are already a reality, particularly in number
theory, and {\em how} they are a reality; such a discussion has to precede
and inform more philosophical considerations.

\subsection{Rigorous computation}

We will be first of all
concerned with the production of numerical results that are
guaranteed to be correct. (The same matter goes by the names of
``reliable computing'', or ``validated numerics''.)
In exploratory work, we may ask, say, for the value of $\sin(0.1)$; however,
a moment's thought shows that the true
numerical value of $\sin(0.1)$ cannot even be stored in a computer --
a decimal or binary expansion has to be cut off at some point. What can
we do?

\subsubsection{Basics on interval arithmetic}

Rigorous computations with real numbers must take into account rounding errors
and other sources of lack of precision. As we were saying, a
generic element of $\mathbb{R}$ cannot even be represented
in a computer. A transcendental number is, in general, given by an infinite
sequence of digits, and computers have only finite space. Furthermore, computers
are built so that they deal with integers, and, by extension, with
rational numbers;
they are fastest when working with rationals whose denominators are powers of
$2$.

{\em Interval arithmetic} provides a way to keep track of rounding errors
  automatically, while providing data types for work in $\mathbb{R}$ and
  $\mathbb{C}$. 
  The basic data type is an {\em interval} $\lbrack a,b\rbrack$, where
  $a,b\in \mathbb{Q}$. For obvious reasons of efficiency, we may
  restrict $a$ and $b$ to be elements of the localization $Q_2$ of $\mathbb{Z}$
  by the powers of $2$, i.e., the ring $Q_2$ of rationals whose
  denominators are powers of $2$. (We also may decide to work with
  $Q_2 \cup \{-\infty,\infty\}$ instead of $Q_2$.)

  A procedure is said to implement a function $f:\mathbb{R}^k\to
  \mathbb{R}$ if, given $B = (\lbrack a_i,b_i\rbrack)_{1\leq i\leq k}$,
  where $a_i,b_i \in Q_2$, the procedure returns an interval $(a,b)$
  containing $f(B)$, with $a,b\in Q_2$. Of course, one would expect a good
  implementation not to return an interval $(a,b)$ much larger than 
  it needs to be.

  Interval arithmetic in this sense was first proposed in the 1950s and 1960s
  (\cite{moore}; see also \cite{Young1931}, \cite{warmus1956},
  \cite{sunaga2009}, among others).
  There are several commonly used open-source implementations.
  (The package ARB \cite{Johansson2017arb} implements a slight variant, {\em ball arithmetic}.)
  A good interval-arithmetic package should implement not just the four
  basic operations and some basic functions such as $\exp$, $\log$ or $\sin$,
  but also as many commonly used transcendental functions as possible.

  The main drawback of interval arithmetic is its time consumption.
  Very roughly speaking, a procedure coded using a good
  interval-arithmetic package
  will typically take about $8$ times as long as the same
  procedure coded without interval arithmetic.

  There is also the fact
  that the underlying floating-point routines must give validated results,
  that is, results with a precise bound on the error term; it is best if
  results are correct up to the last bit, with whatever rounding type
  is specified, as all further bits after the last guaranteed bit are of course
  useless. Correctly rounded floating-point results may feel like a basic
  right, but the fact is that most processors do not guarantee this
  arguable right
  for anything other than the four basic operations,
  and in fact often violate it.
  Thus, transcendental functions have to be implemented in software. See, e.g.,
  \cite{crlibm}. (Of course, it may not be fair to count the need for
  a software implementation as part of the overhead of interval arithmetic;
  we should not use wrong routines implemented in hardware to begin with.)
  
  The speed of computers in our days makes large-scale computations in interval
  arithmetic possible; see, for instance,
  D. Platt's verification of the Riemann Hypothesis for zeroes with imaginary
  part $\leq 1.1\cdot 10^{11}$ \cite{Platt}.

  {\em Alternatives.}  It is possible to avoid interval arithmetic (or rather its usage during the computation, with the attendant
  overhead) by keeping
  track of rounding errors {\em a priori}, that is,
  analyzing carefully to what extent an error of a given size
  in the input to a given procedure can affect the output.
  This is one of the classic subjects of numerical analysis
  \cite{Wilkinson}, \cite{Higham}.
  One clear drawback is that proceeding in this fashion for anything
  outside a small set of well-studied tasks involves saving computer time at the
  expense of human time, and creates one more occasion for human error.

  It would make sense for it to be possible to be able to analyze rounding
  errors a priori {\em with the help of a computer}. (The same goes
  for errors that result from, say, truncating a series.) Systems for doing
  as much exist, and in fact produce formal proofs (q.v.): FPTaylor
  (which produces certificates for the formal system HOL Light),
  Gappa (which gives certificates for the formal system Coq), Real2Float
  (Coq again)\dots

  However, none of these systems can treat loops, at least not if the
  floating-point computation carried out in an iteration of the loop
  depends on the floating-point computations carried out in previous iterations.
  This
  limitation is severe: it makes these systems unusable for one the kinds of computations
  most common in number theory, namely, a computation or verification involving
  sums or some other quantities defined by an iterative or recursive process.
  (For an example, see the discussion of the sum $m(x)$ in \S \ref{sec:asymcomp}.)
  
  On the other hand, systems for analyzing rounding errors can be and have been
  used to verify
  that at least some of the routines used by interval-arithmetic packages are
  correct. Moreover, there are presently tools that, if developed
  further, should become able to give formal certificates for at least some
  kinds of computations with loops in a fairly near future.
  
\subsubsection{Comparisons. Maxima and minima}

In exploratory work, it may be fine to plot two functions
$f,g:I\to \mathbb{R}$, where $I=\lbrack 0,1\rbrack$ (say) and decide that $f(x)<g(x)$ for all $x\in I$
because the graph shows that it is so. Of course,
that would not do as a proof. Fortunately, just as we can let a computer plot
$f$ and $g$, we can let a computer prove that $f(x)<g(x)$ for all $x\in I$.

In fact, it is particularly simple to do so by means of interval arithmetic.
Let $f$ and $g$ be implemented in interval arithmetic. If $f(I)$ and $g(I)$
do not overlap, we obtain either that $f(x)<g(x)$ for all $x\in I$ (if $f(I)$
is to the left of $g(I)$) or that $f(x)>g(x)$ for all $x\in I$ (otherwise).
If $f(I)$ and $g(I)$ do overlap, we divide $I$ in half, and recur:
we do what we have just described
to each half. We stop when we obtain that the statement we are trying to prove
is false on some subinterval, or when we have divided $I$ into subintervals
on each of which the statement is true. While this approach runs into
difficulties for $x$ close to a point $x_0$ such that $f(x_0)=g(x_0)$, we can
usually resolve such a situation by comparing the derivatives or higher
derivatives of $f$ and $g$ at or near $x_0$.

A very similar version of the bisection method in interval arithmetic
(combined, if necessary, with automatic differentiation)
can be used to locate maxima and minima, as well as roots.

\subsubsection{Numerical integration (quadrature)}

There are several well-known ways to estimate an integral and bound the error
in the estimate, starting with, say, the trapezoid rule, or Simpson's rule,
and including Euler-Maclaurin or Gaussian quadrature, for instance.
It is generally possible to implement these methods in interval arithmetic.

Matters can become more complicated if the integrand is not everywhere
differentiable, or even if it not differentiable on the whole {\em closure}
of the interval we are working on. A common case is that of an integrand of the form $|f(x)|$,
where $f$ is everywhere differentiable, but $|f|$ is not differentiable at
all zeroes of $f(x)=0$. That case and several other ones ought to be done
automatically; they still require ad-hoc work at the time of writing.

Even more ad-hoc work is required (for now) if we must compute a complex integral. For instance,
it is possible to show that,
for $P_{r}(s) = \prod_{p\leq r} (1-p^{-s})$ and $R$ the straight path from
$200 i$ to $40000 i$,
\[\frac{1}{\pi i} \int_{R} |P(s+1/2) P(s+1/4)|
\frac{|\zeta(s+1/2)| |\zeta(s+1/4)|}{|s|^2} |d s| = 0.009269 + \text{error},\]
where $|\text{error}|\leq 3\cdot 10^{-6}$. However, ``it is possible'' means
for now ``one can obtain this result by e-mailing ARB's author
(F. Johansson), editing the code he kindly sends you, and then
running it overnight''. It is clear that the situation should (and will)
improve.

\subsection{Combining asymptotics and computations}\label{sec:asymcomp}

Let us now discuss the way in which
the need for computations typically arises in number
theory. (It is not the only way; one can also need, say, verifications of
finite chunks of the Riemann or Generalized Riemann Hypotheses of the kind we
have already mentioned.)
  Often, methods from analysis yield estimates of the following form:
  \begin{equation}\label{eq:experr}
    \textbf{expression}(n) = f(n) + \text{error},\end{equation}
  where $|\text{error}|\leq g(n)$ and $g(n)$ is much smaller than $|f(n)|$
  for $n$ sufficiently large. The question is what to do when $n$ is not
  sufficiently large.

  One answer is, of course, that one may compute $\textbf{expression}(n)$ for
  $n$ not sufficiently large, that is, $n$ smaller than a constant $N$.
  If we reduce our bound $g(n)$, then we are reducing the constant $N$ such that
  $g(n)$ is sufficiently smaller than $|f(n)|$ for $n\geq N$; reducing $N$
  to a reasonable level is thus a challenge for us as analysts.
  Computing $\textbf{expression}(n)$ reasonably efficiently for $n<N$ is
  a challenge for us as programmers, or, what is almost always more
  interesting,  as algorithm designers. Whether the algorithm 
  for computing $\textbf{expression}(n)$ for $n\leq N$ runs in time
  $O(N)$, $O(N^2)$ or $O(N^3)$ is obviously more important than the extent
  to which the code has been optimized.

  An interesting example is given by the proof of the
  ternary Goldbach conjecture, taken as a whole.
  (The first version of the proof is available online, both
  as a series of preprints and as a book draft \cite{Helfbook}; the
  version to be published is in preparation.)
  As is well-known,
  the ternary Goldbach conjecture states that every odd integer $n\geq 7$ can
  be written as the sum of three prime numbers.

  Almost all of the effort involved in the proof went into proving an estimate
  of the form (\ref{eq:experr}) for $\textbf{expression}(n)$ defined
  to be the number of ways one can write $n$ as the sum of three primes
  $p_1$, $p_2$, $p_3$ (counted with certain weights). It was then very
  easy to show that $g(n)<|f(n)|$ for $n\geq 10^{27}$ odd. It followed
  that every odd number $n\geq 10^{27}$ can be written as the
  sum of three primes, i.e., ternary Goldbach holds for $n\geq 10^{27}$.
  It remained to show that every odd number $7\leq n\leq 10^{27}$ can also be
  written as the sum of three primes; then the ternary Goldbach conjecture
  would follow.

  Checking each odd number $\leq 10^{27}$ is out of current computational reach.
  However, there is the following well-known trick. We can easily construct
  an increasing succession of primes, all of them $<N$ but for the last
  one, such that the difference between any two consecutive primes in the list
  is at least $4$ and at
  most $M-2$, say; the time taken is not much larger than $O(N/M)$.
  (Platt and I checked that for $N = 8.875\cdot 10^{30}$ and
  $M = 4\cdot 10^{18}$; by now, there is even a formal proof of the same
  for $N=10^{28}$ and $M=4\cdot 10^{18}$ (see \cite{ThGr}).) Then we can use a
  verification by brute force
  that the strong Goldbach conjecture holds for all even
  numbers $4\leq n\leq M$, i.e., every even number $4\leq n\leq M$ can be
  written as the sum of two primes. It then follows that every odd number
  $7\leq n\leq N$ can be written as the sum of three primes: just let $p$
  be the largest prime in the sequence such that $n-p\geq 4$, and apply
  strong Goldbach to $n-p$, which has to be at least $M$. We obtain that
  $n-p=p_1 + p_2$.
  Thus $n = p + p_1 + p_2$, and so we are done. The total time taken is
  not much larger than $O(M+N/M)$, so we can simply set $M=\sqrt{N}$,
  and obtain an algorithm running in time not much more than $O(\sqrt{N})$,
  which is not too much for $N=10^{27}$.
  
  As it happens, it had already been checked that strong Goldbach holds for
  all even $n\leq M = 4\cdot 10^{18}+2$ (\cite{e2014empirical},
  plus a trivial check for
  $4\cdot 10^{18}+2$). Hence, we can simply set $M=4\cdot 10^{18}+2$, and the proof
  of the ternary Goldbach conjecture is done.

  Let us see a different, much smaller example. We would like to bound
  the sum $m(x) =\sum_{n\leq x} \mu(n)/n$ for general $x>0$. (That sum and its
  variants play an important role in the proof of ternary Goldbach.)
  The strongest
  explicit bound is due to Ramar\'e \cite{ramare2015explicit}, who proved that
  \begin{equation}\label{eq:ramar}
    |m(x)|\leq \frac{0.0144}{\log x}\end{equation}
  for $x\geq 96955$. (We would have a bound of $O_\epsilon(1/x^{1/2-\epsilon})$ if
  we assumed the Riemann hypothesis, but that is out of reach; there
  are also unconditional non-explicit bounds that are qualitatively
  much stronger than (\ref{eq:ramar}) as $x\to \infty$, but there are
  serious obstacles to making them explicit with constants reasonable enough
  for the bounds to be usable.)

  Clearly, we can compute $m(x)$ for all small integer values of $x$.
  I wrote a C program using interval arithmetic to do so, and obtained, after
  a couple of weeks on a rather ordinary PC, that
  \begin{equation}\label{eq:smallmx}
    |m(x)|\leq \sqrt{2/x}
  \end{equation}
  for all real $0<x\leq N = 10^{14}$, with $0.569449$ instead of $\sqrt{2}$ if
  $x\geq 3$ is assumed. Note that, if a standard conjecture holds, the bound
  (\ref{eq:smallmx}) is actually too strong to hold for all $x$.

  Is it overshooting to test (\ref{eq:smallmx}) for all $x\leq 10^{14}$? Perhaps, but, since the
  obvious algorithm for checking (\ref{eq:smallmx}) for all $x\leq N$ takes
  time essentially linear on $N$, we might as well allow ourselves a large $N$.
  Also to the point -- it is possible to combine (\ref{eq:ramar}) and
  (\ref{eq:smallmx}) to prove a result valid for all $x$, small and large:
  given any $y>1$ and any $0<x\leq y$,
 \[|m(x)|\leq \sqrt{\frac{2}{x}} + 0.0144\cdot \frac{y^c}{\log y} \frac{1}{x^c},\]
  where $c=1/\log N= 1/\log 10^{14}$.
  The bound here is a sum of powers of $x$, something highly convenient for
  several purposes. The larger $N$ is, the better the bound.

  Of course, in this example, it was easy to give a fast algorithm. For
  some other sums I had to deal with, the obvious algorithm would have run in
  time $O(N^2)$ or $O(N^3)$, and finding an algorithm that would run in,
  say, time $O(N \log N)$ and space roughly $O(\sqrt{N})$ was a challenge.

\subsection{Error avoidance}

In practice, computer errors in the strict sense
are barely an issue, except perhaps for very large computations. (By
``computer error'' we mean a malfunction in a correctly designed circuit.
There is also the issue of incorrectly designed circuits, such as those in the
original Pentium chip, whose flaw, as it happens, was found through a
computation in number theory (\cite{Pentium}, \cite{cipra1995number}).
Such situations are nowadays largely avoided
through {\em formal verification} -- another interesting area.)
Computers have long had various inbuilt tools (e.g., checksums) for error
detection and correction. Moreover, there is the obvious fact that the same
computation can be run time and again, possibly on different hardware.
Such a thing would be cumbersome or expensive
only if the computation were very large indeed.
(Admittedly, that is also the case when the probability of computer
error might not be overwhelmingly negligible.)

Almost all of the time, the actual issue is human error: errors in programming,
and errors in input/output. Of course, human beings
also make mistakes when they don't use computers. There are conceptual
mistakes, and then there are silly errors,
especially in computations or tedious case-work. Presumably, whatever discipline
we adopt to avoid errors in programming and input/output can also be adapted
to weed out silly mistakes in general.

Before we discuss formal proofs, we should touch upon the more humble matter
of everyday discipline.\footnote{Specialists in formal proofs remind me that
  ``everyday discipline'' in the sense I am about to discuss is also an issue
  for them.}
Here I can simply describe my current practice,
particularly in the context of my proof of the ternary Goldbach
conjecture. 

When facing computational tasks, I tend to program a great deal at first.
However, I
try to minimize the number of programs used in the submitted version of
a book or paper, and keep them short. Nowadays I classify tasks in two
categories. For time- and
space-intensive computations, I am whittling down my code to a few programs
in C, submitted or to be submitted to referees. The programs are or will be
available upon request. 

For small computations, I now use Sage/Python code, almost all of it
included in the TeX source of the book. (The code for small computations that
take more than a couple of seconds has been submitted together with the programs
in C just mentioned.) It is an easy matter to keep Python
code brief and readable. Thanks to SageTeX, input and output are automated,
in that the output of Sage/Python code is displayed automatically in the 
file that TeX outputs. For instance, given two variables called \texttt{result1}
and \texttt{result2}, with values $5$ and $7$, we can type the TeX code
\begin{verbatim}
Hence,  $f(x)\leq \sage{result1 + result2}$.
\end{verbatim}
and obtain
\begin{quote}
  Hence, $f(x)\leq 12$.
\end{quote}
In this way, the possibility of human error in copying input/output and
updating material is reduced. All code can be run afresh between two
compilations of the TeX file; indeed, I run it afresh (and have to run it
afresh) after any change to any code in the file. Someone downloading the
source code may decide to run the code or not; it is very easy to do so
for anyone with Sage installed. It is also very easy to display all code;
all that is involved is switching
a single flag at the beginning of the source code.

There is a question that cannot be avoided, namely, whether
it is right to use a large computer-algebra system in a proof. It is in
general unsatisfactory to say ``By M.,\dots''\footnote{Here
  M. stands for Mammoth, or anything else that is closed-source and very large.} before a claim, as some
papers do. Consequently, in the first version of my proof of ternary
Goldbach, I used Sage only in an exploratory role, and did all coding in C.
At the same time, Sage is open-source, peer-reviewed and highly modular
(though some of its larger, older components may not themselves be highly
modular or peer-reviewed in quite the same way).
Because Sage is modular, we depend only on the correctness of the relatively
small parts of it that we use in the proof; in my case, that amounts to
Python, basic symbolic algebra and ARB.
In part for this reason, and in part for the reasons outlined above
(code readability, code organization,
reduction of human error in input/output),
I decided that the advantages of using Sage/Python for small computations
that can be embedded in the TeX code overwhelmed any reasonable misgivings
one might have about using Sage in a proof in this particular way.

\subsection{Automated proofs. Formal proofs}

\subsubsection{Automated proofs in practice}

Informally speaking,
G\"odel's incompleteness theorem states that, in a general axiomatic
system -- in particular, in any finite axiomatic system that includes basic
arithmetic, including exponentiation -- there are truths that cannot be proved
starting from the axioms. In particular, we cannot ask a machine to prove
all true statements.

Moreover, even in axiomatic systems that {\em are} complete, the problem
of finding a proof may be (and in fact is) computationally hard. Indeed,
the first problem to be proven to be NP-complete was that of Boolean
satisfiability (SAT); in other words, if we could determine rapidly whether the
value of the variables in a Boolean formula can be set so that the formula
is true, we could solve rapidly
{\em any} problem in a very broad class (essentially, any problem whose
solution can be verified rapidly).

Does that mean that fully automated theorem proving is hopeless? Not
necessarily. First, a machine could prove {\em some} true statements,
even if it cannot prove every true statement; after all, human beings are in
that same situation. Second, a computer can proceed by heuristics to solve
quickly some cases of a problem that is computationally hard in general.
Third, a problem may be small enough (in terms of its number of variables, say)
that an inefficient algorithm running on a computer can still solve the
problem in a
reasonable time, even when its solution would be non-obvious or cumbersome to a
human.

I have had little experience with automated theorem proving -- and in fact
was surprised to see that it could be used in practice at all. At some point,
I wanted to prove the following lemma: 
 for all  $0<x\leq y_1,y_2<1$ with $y_1^2\leq x$, $y_2^2\leq x$,
\begin{equation}\label{eq:odmalicka}
1 + \frac{y_1 y_2}{(1-y_1+x) (1-y_2+x)} \leq
 \frac{(1-x^3)^2 (1-x^4)}{(1- y_1 y_2) (1 - y_1 y_2^2) (1 - y_1^2 y_2)} .
\end{equation}
Now, this is a statement in the theory of real closed fields, which is in fact
complete; a proof has to exist. It turns out that there is a freely available
program QEPCAD \cite{QEPCAD} implementing an algorithm that proves statements
in the theory of real closed fields (by {\em CAD}, that is,
cylindrical algebraic
  decomposition \cite{MR0403962}). Can it deal with (\ref{eq:odmalicka})?
Not on my desktop, and that is not surprising: the computational
complexity of CAD is doubly exponential in $n$, with base proportional
to the maximal degree of the polynomials involved. However,
thinking a little, it is not hard to eliminate a variable in
(\ref{eq:odmalicka}), in that one can show that the maximum of the
left side minus the right side has to be attained when either
$y_1 = y_2$, $y_i = \sqrt{x}$ or $y_i = x$ holds for at least one of $i=1,2$.
QEPCAD proves the resulting inequality in two variables with ease.

Since (\ref{eq:odmalicka}) has quantifiers of only one kind
($\forall x \forall y_1 \forall y_2$), there are alternative algorithms
that solve the problem in time exponential on the number of variables
(see \cite[Ch. 11]{opac-b1124307} and references therein).
However, there does not seem to be a practical, reasonably efficient
implementation of the exponential-time algorithm commonly available just yet.

In the end, I found an alternative way to prove what I truly wanted, and so
I no longer needed (\ref{eq:odmalicka}) at all. It was still interesting to
learn that automatic provers could sometimes establish useful auxiliary lemmas. (Of course, in some sense,
the simple combinations of the bisection method and interval arithmetic
we have already discussed also fall into this category.) It was also
interesting to learn of the extent to which
practice still fell behind theory -- and
of how ``computationally hard'' does not always mean ``hopeless'' in practice.

Again, my personal comments should be understood as coming from the possibly
na\"ive perspective of a number theorist. I have recently learned that
Boolean-satisfiability (SAT) solvers are being used intensively in
combinatorics. There, one of the main issues is error in the human/computer
interface: the interface of an automated solver is not necessarily
intuitive, and some claimed proofs in the literature are the result of
inputting the problem into a SAT solver incorrectly.

At any rate, it seems clear that the role of automated provers in the
foreseeable future will be one of {\em assistance} to the human prover.
Thus we are led to our next, broader subject.

\subsubsection{Formal proofs and proof assistants}

Contrary to what we may sometimes tell ourselves,
what we call a proof in our everyday practice is not quite the same as what
we would call a proof in a logic course. The latter -- called a
{\em formal proof} for clarity -- is a sequence of symbols whose
correctness is a purely
syntactic property that can be checked by a monkey grinding
an organ. In contrast, a proof, for the working mathematician, is a convincing
argument that can in principle be turned into a formal proof.

Until relatively
recently, ``in principle'' came with an enormous caveat: none but the
simplest sort of proofs was turned into formal proofs. Our imaginary
organ and monkey were first replaced by a (real) computer system in 1967
(de Bruijn's AUTOMATH \cite{de1970mathematical}),
but the task of producing a formal proof remained solely
in the hands of humans, and was very cumbersome.

The situation has
gradually changed thanks to {\em proof assistants}, programs whose role it is
to help a human being write a formal proof. There is now a number of such
assistants: Mizar, Coq, Isabelle, HOL Light,\dots

One notable recent success was the second proof of Hales'
sphere packing theorem. As is well-known, the first proof
(\cite{Hales2005}; see also \cite{lagarias2011kepler}) 
was computer-assisted in the more traditional ways discussed in previous
sections. That first proof
met with some misgivings: it involved a great deal of case-work
(thus increasing the possibility of human error, in part because
the
refereeing process became very onerous), and the computer code 
ran up to about 40000 lines (making the code practically impossible to referee).
The second-generation proof \cite{Halesetal2017a} is a formal proof,
written and verified by means of the Isabelle and HOL Light proof assistants.

In which fields of mathematics is such an effort now feasible?
Large areas, including much of basic real analysis, remained uncovered until
recently. As far as analytic number theory is concerned: the first formal
proof of the Prime Number Theorem \cite{avigad2007formally} was constructed
in 2005, by means of Isabelle; it is based on Selberg's elementary proof,
which is
often seen as more difficult, or less natural, than more traditional proofs
by complex analysis. The first analytic proof, based on Newman's simplified
version of the traditional approach \cite{newman1980simple}, was given in
\cite{harrison-pnt}; it relied on the fact that some of real and complex analysis had become available in HOL Light \cite{harrison2005hol}, \cite{harrison2007formalizing}, and also necessitated giving
formal proofs of some basic facts about the Riemann zeta function. 

In analytic number theory, we are, then, perhaps just past the beginning.
The next natural step would be to formalize much of a first-year graduate
textbook -- say, the results in \cite{MR0217022} and \cite{MR2378655},
possibly with different
proofs, together with some sieve theory, and also \cite{zbMATH03968684};
then we would need large parts of \cite{MR2061214}.
(Replace chapters of the older books here by more modern sources when needed.) Then we would be able to
start working on newer or new material. Of course, one can also proceed
backwards, setting oneself a challenge (several of the major results in
\cite{MR2061214} would do nicely) and working backwards from it,
proving whatever basic results one needs, much as
Gonthier et al. did with the Feit-Thompson Theorem \cite{Gonthieretal2013a}.

This last strategy is, in a sense, similar to what I did in the end for
ternary Goldbach; I had to ask about, learn about, and prove explicit
analogues of many basic results in analytic number theory. Sometimes, to my
na\"ive surprise, I had to do without a standard technique or result, since
no practical explicit analogue existed or could realistically be proved.
We will presumably face the same kind of challenge when we try to give 
formal proofs of the main results of twentieth-century number theory,
including, why not, Vinogradov's three-prime theorem.

What about formal proofs of explicit results, or of statements whose proofs
make crucial use of explicit intermediate results? Why not a formal proof of
ternary Goldbach? Why not indeed, in the long run, but we have to be realistic
about the fact that both making a proof explicit and making it formal take
plenty of work, and of course can also lengthen the proof greatly. It may be
that the de Bruijn factor -- that is, the informal quantity defined as the
ratio of the length of the formal version of the proof to the original, conventional (``informal''?) version, however stored -- is lower for some explicit results, as their proofs tend to
include a level of detail not always needed for non-explicit results. Time will
tell. Of course one can also say that formal proofs are more sorely needed
for explicit results; while, in general,
they do not elicit the same suspicion as results with plenty of case-work,
they may be more likely to be in the ``fixable, but incorrect as stated''
category than their non-explicit counterparts. (Asymptotic notation is
a carpet under which both known and unnoticed dirt can be conveniently swept.)
Thus it would be a very worthwhile task -- for the fairly near future --
to give formal proofs for basic explicit
estimates that are used again and again.

\subsection{Final considerations}

A point often made in connection with computer-assisted proofs is that the
purpose of proof is not just to establish truth, but to demonstrate and
advance understanding. The worth of a proof understandable to no one
would thus be limited.

This viewpoint is valid, but may not really be specially relevant to
computer-assisted proofs as they are currently developing. It could certainly
be a point against automated provers, particularly if their output is not
intelligible. However, in most fields,
automated provers seem likely to continue playing at most an auxiliary role, 
proving small lemmas that would be cumbersome and not particularly
enlightening for a human to prove. (Does (\ref{eq:odmalicka}) have any
``meaning''?)
Computer-assisted formal proofs are a
different matter: there, we start with a proof in the everyday sense, that is,
a proof that is produced by human beings and understandable to a (hopefully
proper) superset of the same; we then formalize it, with the help of a computer,
having the more solid establishment of truth as the primary or sole aim.

As for case-work: while it is often locally easy to understand (so to
speak), it is true that it
can feel meaningless, unaesthetic, and an invitation to error by
means of tedium. However, that is an issue with case-work in general, and not
with computer-based case-work in particular. Computers can, at least, play
a role in reducing or eliminating error from case-work.

Moreover, the kind of finite verification typical of number theory
cannot really be said to be case-work in this sense. Instead of having many
slightly different cases, we have typically a single equality or inequality
that must be verified for very many integers, or for all values of a variable
under a certain threshold. Moreover, what we verify is often far from meaningless: if a computation verifies that the first $10^9$ non-trivial zeroes of the
Riemann zeta function lie on the critical line, or that $\sum_{n\leq x}
\mu(n)$ is bounded by not much more than $\sqrt{x}$ for all $x\leq N$,
we are verifying finite parts or finite consequences of standard conjectures
that we have very good reasons to believe in -- and good reasons
to believe out of reach.

The assistance of computers can lead us both to results that were previously
unattainable and to a higher standard in certainty and rigor.
It may be some time before the standard of rigor set by formal proofs
becomes widespread, but we have come to the point where computers, correctly
used, can quell misgivings rather than give rise to them. 


{\bf Acknowledgements.} Many thanks are due to Guillaume Melquiond, Assia Mahboubi and Victor Magron for their feedback. Funding from his Humboldt professorship and from his ERC Consolidator grant 648329 (GRANT) is gratefully acknowledged.

\section{Machine-Assisted Proofs in Group Theory and Representation Theory:
Pham Huu Tiep\protect\footnote{The author gratefully acknowledges the support of the NSF (grants DMS-1839351 and DMS-1840702). He also thanks 
Gabriel Navarro and Eamonn O'Brien for helpful comments on the topic of this discussion.}}



Typically, many proofs in mathematics rely on {\it mathematical induction}. In group theory and representation theory,
this inductive approach often follows a modified strategy, which can be described as follows. Suppose the goal is to prove a 
certain statement $\PC$ concerning a (finite or algebraic) group $G$. 

\begin{enumerate}[\rm(i)]
\item Then the first step is to prove a {\it reduction theorem}
to reduce to the case where $G$ is (very close to be) {\it simple}, using perhaps the
{\it Classification of Finite Simple Groups}. (One should note that, these reduction theorems usually require one to 
establish a much stronger condition $({\mathcal P}^*)$ for simple groups than the original condition $\PC$. 
See \cite{IMN} and \cite{NT} for some recent reduction theorems.)

\item Next, one works out a uniform proof, which handles the simple groups $G$ of {\it large enough} order.

\item The {\it induction base} is then to treat all ``small'' simple groups $G$.
\end{enumerate}

In either strategy, the induction base usually needs a completely different treatment, rather than the uniform case of large groups,
which often involves the, sometimes massive, use of computer calculations.
 
Let us illustrate this on the example of the proof of the {\it Ore conjecture} \cite{O}:

\begin{conj}[\textbf{Ore, 1951}]\label{ore}
Every element $g$ in any finite non-abelian simple group $G$ is a commutator,
i.e. can be written as $g = xyx^{-1}y^{-1}$ for some $x,y \in G$.
\end{conj}

Many important but partial results on the conjecture were established by Ore himself, Miller, R. C. Thompson, Neub\"user-Pahlings-Cleuvers, and most notably, Ellers-Gordeev. The conjecture was finally proved in \cite{LBST}:

\begin{thm}[\textbf{Liebeck-O'Brien-Shalev-T, 2010}]
Conjecture \ref{ore} holds for all finite non-abelian simple groups.
\end{thm}

Even building on all previous results, the proof of this ``LOST-theorem''
is still 70-pages long. So how does this proof go? A detailed account of it was given in the 2013 Bourbaki seminar \cite{M}.
A key ingredient of the proof is the following formula, where $\Irr(G)$ denotes the set of all complex irreducible characters of 
the finite group $G$:

\begin{lem}[\textbf {Frobenius character sum formula}]\label{frob}
Given a finite group $G$ and an element $g \in G$, the number of pairs $(x,y) \in G \times G$ such 
that $g = xyx^{-1}y^{-1}$ is 
$$|G| \cdot \sum_{\chi \in \Irr(G)}\frac{\chi(g)}{\chi(1)}.$$
\end{lem}

So in order to show that a given element $g \in G$ is a commutator, one just needs to show that 
$\sum_{\chi \in \Irr(G)}\chi(g)/\chi(1) \neq 0$. Now, let $G$ be one of the groups in the induction base for
the proof \cite{LBST} of the Ore conjecture. 

\begin{enumerate}[\rm(a)]
\item In many cases, the character table of $G$ is well known, either 
published and/or publicly available, in which case we can just use Lemma \ref{frob}. 

\item In the remaining case, where the character table of $G$ is not available, but if $|G|$ is not too large,
then we construct the character table of $G$ and then proceed as before. To construct the character table of
such a group $G$, one starts with a ``nice'' presentation or representation of $G$. Then one can try to use various operations with
group characters to produce enough characters of $G$ to generate the full group $\ZZ\Irr(G)$ of complex characters of 
$G$. With respect to the usual inner product
$$[\alpha,\beta] = \frac{1}{|G|}\sum_{x \in G}\alpha(x)\overline{\beta(x)},$$
$\ZZ\Irr(G)$ is a Euclidean lattice, whose minimal vectors are (up to sign) precisely the {\it irreducible} characters of $G$ that we are after.  
So one can try to follow, say the {\it LLL-algorithm}, to find these irreducible characters. In practice (and in \cite{LBST}), one can use 
Unger's algorithm \cite{U}, implemented in {\sf MAGMA} \cite{Magma}.

\item But some of the groups $G$ in the induction base of \cite{LBST}, like ${\mathrm {Sp}}_{10}(\FF_3)$, $\Omega_{11}(\FF_3)$, or $\mathrm{U}_6(\FF_7)$, are still too big for the computation in (b). For these too-big groups $G$, we implement another 
strategy. Namely, for any given $g \in G$, we  run a randomized search for $y \in G$ such that $y$ and $gy$ are 
conjugate. (In fact, one needs to do it only for one representative of each conjugacy class of $G$. So, for the largest 
sporadic simple group, the {\it Monster}, of order about $8 \cdot 10^{53}$, one needs to work with only $194$ such representatives.)
Once such a $y$ can be found, then we have $gy = xyx^{-1}$ for some $x \in G$, i.e. $g = [x,y]$, which is what we wanted to show! 
\end{enumerate}

The first obvious question that arises is: How {\it long} was this computation? All in all, it took us about 150 weeks of CPU time of a 2.3GHz computer with 250GB of RAM to complete (by April 2008) all the computations needed for the proof in \cite{LBST} of the Ore conjecture. (Certainly, this amount of CPU time could be significantly reduced with the better computational algorithms available now, ten years later.)

The second, and more important, question is: How {\it reliable} is this computation? In the cases of (a) and (b) where we used or computed the character table of $G$,  the relevant character tables have been subjected to various checks which attest to their accuracy.  The tables are also publicly available in character table libraries, so they can be checked by others and used for independent verification of the conjecture. 
Next, in the larger cases of (c), the randomized computation was used to find $y$ that $gy$ and $y$ {\it should} be $G$-conjugate for
a given $g$, and then one checks directly (using the given presentation or representation of $G$) that $gy$ and $y$ are {\it indeed}
$G$-conjugate. One should notice that this is quite different in nature to other machine-assisted proofs which reduce an elaborate proof to many cases -- each is then decided by machine, often reporting ``yes'' or ``no'' to the existence of some object.

\bigskip
In group theory and representation theory, as in many other areas of mathematics, perhaps even more important than 
the machine-assisted proofs are

\begin{itemize}
\item {\it machine-assisted discovered theorems}, and  
\item {\it machine-assisted discovered counterexamples}.
\end{itemize}
Let us mention a couple of examples of these two kinds.

\begin{itemize}
\item The {\it Galois-McKay conjecture} 
was formulated by Navarro in \cite{N} after many, many days of computing in {\sf GAP} \cite{GAP}.

\item Also after some long experimental computations with the symmetric groups $\SSS_n$, $n \leq 50$, Isaacs, Navarro, and I found a natural 
{\it McKay correspondence} (for the prime 2), which should hold for all symmetric groups, and which was subsequently 
proved in our joint paper with Olsson \cite{INOT}.

\item An old conjecture in Character Theory states that if a finite group $G$ is rational (that is, all $\chi \in \Irr(G)$ are 
rational-valued), then so are the Sylow $2$-subgroups of $G$. However, after some long (but directed!) search using 
the ``SmallGroups'' databases contained in the computer packages \cite{GAP} and \cite{Magma}, 
Isaacs and Navarro have been able to find two counterexamples of order $2^9 \cdot 3$ to this conjecture, see \cite{IN}. 
\end{itemize}

\section{Lu\'\i{}s Cruz-Filipe}
The awareness of the mathematical community towards the use of computers in
proofs of mathematical results started in 1976 with the announcement of the
proof of the Four-Color Theorem~\cite{AppelHaken1976a}.

However, there are earlier examples of proofs that were partially or completely
done by a machine. An interesting example is Floyd and Knuth's proof of
optimality of the 16-comparator sorting network on 7 inputs (Theorem 5
in~\cite{FloydKnuth1973}), which starts:
\begin{quote}
This theorem was proved by exhaustive enumeration on a CDC G-21 computer at
Carnegie Institute of Technology in 1966.
\end{quote}

Optimality of sorting networks is an interesting problem in combinatorics. The
question we want to answer is: what is the length $S(n)$ of the shortest
sequence of compare-and-swap operations (i.e., atomic operations consisting of
sorting a pair of values) that will sort all inputs of a given length $n$?
Floyd and Knuth's work cited above addresses this problem for $2\leq n\leq 7$.
Except for the values of $S(4)$ and $S(6)$, which are derived from those of
$S(3)$ and $S(5)$ by application of a more general theorem, all values are
derived by exhaustively enumerating the possible sequences of length $S(n)-1$
and showing that, for each of them, there is an input of length $n$ that they do
not sort. As the authors note in their conclusion, this method is however
``quite unsatisfactory for higher values of $n$'': the number of cases that need
to be analyzed for establishing $S(7)$ is unmanageable for a human being, and
lay at the limits of what could be computed in 1966.

\paragraph{Two styles of proving.}

Floyd and Knuth's proof, as well as the proof of the Four-Color Theorem, are
examples of one family of machine proofs: they establish a property by means of
an \emph{ad-hoc} computer program that is written for that specific purpose.
Nowadays, this approach is not very common: verifying such proofs also demands
that the program be verified, a task that does not fall under the usual scope of
the peer-reviewing process. Instead, it is generally considered preferable to
use a \emph{theorem prover}: a general-purpose program that can construct and/or verify proofs in a particular logic.

Trusting a proof produced with the help of a theorem prover also has some
implications. First, one still needs to trust a computer program -- the theorem
prover itself. The difference with respect to using \emph{ad-hoc} programs is
that we are now considering a program that has been subject to much wider
scrutiny, and it can be reasonably argued that the ensuing proof is as
trustworthy as a published mathematical proof that has been subject to
peer-reviewing. The other point that needs to be checked is that the encoding of
the actual mathematical problem in the theorem prover's logic is correct: this
is generally accepted as part of the peer-reviewing process, as these encodings
are typically included and discussed in submitted articles. (This aspect is also
an issue when using \emph{ad-hoc} computer programs in proofs, but in that
scenario the encodings used tend to be much more direct.)

Theorem provers come in a wide variety of styles and flavors, ranging from very
general purpose to more specifically tailored for a particular family of
problems. They use different logics, and are often uncomparable in their
usefulness. A very encompassing overview of the world of theorem provers can be
found in~\cite{Wiedijk2006}.

\paragraph{Machine-assisted proofs today.}

Forty years after the announced proof of the Four-Color Theorem,
machine-assisted proofs are everywhere. Arguably, their widest area of
application is hardware and software verification, of which Floyd and Knuth's
problem is an example. In an area where we are increasingly dependent on
computers to perform critical tasks -- from controlling air traffic to
administering medicine to patients -- it is more and more important that there
be no errors in the execution of those tasks, whether due to programming errors
or to hardware flaws.

As regards hardware, many useful properties can be verified in a fully automated
way~\cite{Kropf1999}. As such, formal verification is becoming more commonplace
in industry. Software verification tends to be more complex, but the number of
success stories is increasing. In recent years there have been large projects
addressing formalization of large fragments of widely-used programming
languages, making it possible for non-experts to experiment with proving
properties of the programs they write.

Much more challenging is the formalization of mathematical proofs in a
computer. These proofs tend to be much less mechanic and systematic, requiring
constant interaction between the computer and an expert -- where the expert
``guides'' the computer through the main steps of the proof. Mathematics also
tends to build upon itself, and researchers often find that the biggest
challenge in formally verifying an ``interesting'' result lies on the
unavailability of libraries formalizing the relevant underlying theories.
Nevertheless, recent formalizations of e.g.~the proof of the Four-Color
Theorem~\cite{Gonthier2008} and the proof of Kepler's
conjecture~\cite{Hales2014a} show that this goal is not hopeless. (A recent
contribution to this list is a formal proof of Floyd and Knuth's results on
sorting networks~\cite{lcfEtAl2017a}, along with the recently established value
of $S(9)$~\cite{CodishEtAl2016}.)

A different approach to machine-verified proofs of mathematical results, which
some may argue is less elegant, capitalizes on the sucess of SAT solvers. A SAT
solver is a program that tackles the Boolean satisfiability problem: given a
propositional formula, decide whether it has a satisfying assignment. Although
this problem is NP-complete, years of intensive investiment in researching
efficient techniques for solving it has resulted in extremely powerful tools
that can solve formulas with millions of variables and billions of clauses. To
test the limits of these programs, several researchers have experimented with
encoding open mathematical problems (typically from combinatorics) as
propositional formulas, which were then successfully proven to be
unsatisfiable~\cite{HeuleKullmann2017}.

Until recently, SAT solvers were viewed by some with skepticism: their
complexity made them impossible to analyze in practice, and in the case a
formula was claimed to be unsatisfiable no independently verifiable guarantee of
this fact was provided. The situation changed recently, and the majority of
today's SAT solvers also output traces that allow independent, formally verified
systems to check proofs of unsatisfiability of propositional
formulas~\cite{lcfEtAl2017b,lcfEtAl2017c}.

\section{Questions}
\begin{description}
\item[Q]Isn't there a problem with proprietary systems (e.g. Magma, Maple or Mathematics)? You might not be able to find referees who have them.
\item[A]There is a problem here, but it is more about rarity than about the proprietary nature, and in fact you are more likely to find a referee who knows (and has) one of these than you are some rare open-source package. The real problem is the proprietary nature of the algorithms. ``Yes, this large piece of code, which may be proprietary or open-source, gives me this result, but do I trust it?''  For computer algebra, the question is discussed in \cite{Davenport2017z}.
\par On the other hand, it is also the case that Sage has made installing open-source mathematics programs easier, in so far as it seems to incorporate most of them. 
\item[Q]What about a system that produces verified code, which is then compiled and run?
\item[A]That's a good question, and there is some of this in Flyspeck \cite{Halesetal2017a}.  See also Paulson's MetiTarski project \cite{AkbarpourPaulson2010}.
Of course, one would need a verified compiler, but such things exist these days.
\item[Q]Is there much consistency between journals on how these proofs are treated?
\item[A]There isn't even much consistency within a given journal. \cite{Hales2005} and \cite{Maynard2015a} were both in \emph{Annals of Mathematics}, yet seem to have been treated differently. There is no caution on \cite{Maynard2015a}, and it was published much more rapidly than  \cite{Hales2005}. Conversely the computer programs underpinning \cite{Hales2005} are on the \emph{Annals of Mathematics} website, whereas \cite{Maynard2015a} simply says ``An ancillary Mathematica(R) file detailing these computations is available alongside this
paper at \url{ www.arxiv.org}''. The computations in \cite{Maynard2015a} were relatively simple, and carried out exactly: what would \emph{Annals of Mathematics} have made of \cite{Maynard2013a}?
\end{description}

\bibliography{ICMPanel2}
\end{document}